%

\input amstex
\magnification=1200
\loadmsam
\loadmsbm
\loadeufm
\loadeusm
\UseAMSsymbols

\hsize=6.9truein
\hoffset=-0.11truein
\vsize=8.9truein
\voffset=-0.2truein

\def\leftitem#1{\item{\hbox to\parindent{\enspace#1\hfill}}}

\def\boxit#1#2{\hbox{\vrule
	\vtop{%
	\vbox{\hrule\kern#1%
	\hbox{\kern#1#2\kern#1}}%
	\kern#1\hrule}%
	\vrule}}

\def\leaderfill{\leaders\hbox to 1em{\hss.\hss}\hfill}

\parskip=\medskipamount
\document

\input epsf


\centerline{\bf Products of Positive Dehn-twists on Surfaces}

\centerline{ Feng Luo}

This is an appendix to a paper by Michael Freedman. 
The purpose of this note is to prove the following results.

{\bf Theorem 1.} \it Suppose $\Sigma_g$ is a closed surface with a hyperbolic
metric of injectivity radius $r$. There exists a computable
constant $C(g,r)$ so that
each isometry of $\Sigma_g$ is isotopic to a composition of positive
and negative Dehn-twists $D^{\pm 1}_{c_1} ....  D^{\pm 1}_{c_k}$
where $k \leq C(g,r)$ and the length $l(c_i)$ of $c_i$ is  at
most $C(g,r)$ for each $i$. \rm


Call  a self-homeomorphism of the surface \it positive \rm  if
it is isotopic to a composition of positive Dehn-twists.

{\bf Theorem 2.} \it Suppose $\Sigma_{g,n}$ is a compact orientable surface
of genus $g$ with $n$ boundary components. Let
\{$a_1, ...., a_{3g-3 + 2n}\}$ be a 3-holed sphere decomposition
of the surface where $\partial \Sigma_{g,n} = a_{3g-2+n} \cup ... \cup
a_{3g-3+2n}$. Then each orientation preserving homeomorphism of the
surface which is the identity map on $\partial \Sigma_{g,n}$ is
isotopic to a composition $qp$ where $p$ is positive and $q$ is
a composition of negative Dehn-twists on $a_i$'s. \rm

The basic idea of the proof of theorem 1 suggested by M. Freedman
is as follows. Let $f$ be an isometry of the surface. Choose a
surface  filling system of simple geodesics $\{s_1, ..., s_k\}$
whose lengths are bounded (in terms of $r$ and $g$). Since
the lengths of $s_i$ and $f(s_j)$ are bounded, the intersection
numbers between any two members of \{$
s_1, ..., s_k, f(s_1), ..., f(s_k)\}$ are bounded. Now the
proof of Lickorish's theorem in [Li] is constructive  and depends only on
the intersection numbers between simple loops. 
Thus, one produces a bounded number of simple loops of bounded lengths
so that the composition of   positive or negative 
Dehn-twists on them sends $s_i$ to $f(s_i)$.
This shows that $f$ is isotopic to the composition.

The proof below follows the Freedman's sketch. We shall choose the
surface filling system to be of the form $\{
a_1, ..., a_{3g-3}, b_1, ..., b_{3g-3}\}$ where $\{a_i\}$ forms
a 3-holed sphere decomposition of the surface so that $l(a_i)
\leq 26(g-1)$ (Bers' theorem) and $b_i$'s have bounded lengths
so that $b_i \cap a_j = \emptyset$ for $j \neq i$.
Then we establish a controlled version of Lickorish's lemma
(lemma 2 in [Li]) by estimating the lengths of loops involved
in the Dehn-twists. 

We shall use the following notations and conventions. Surfaces are
oriented. If $a$ is a 
simple loop on a surface, $D_a$ denotes the positive Dehn-twist along
$a$ and $l(a)$ denotes the length of the geodesic isotopic to $a$.
Two isotopic simple loops $a$ and $b$ will be denoted by $a \cong b$.
Given two simple loops $a,b$, \it their geometric intersection number \rm
denoted by $I(a,b)$ is min$\{ |a' \cap b'| | a' \cong a, b' \cong b\}$.
It is well known that  if $a, b$ are two distinct simple geodesics,
then $|a \cap b| = I(a, b)$. We use $|a \cap b | = 2_0$ to denote
two simple loops $a, b$ so that $I(a, b) = |a \cap b| =2$ and their
algebraic intersection number is zero.

To prove theorem 1, we begin with the following.

{\bf Proposition 3.} \it Suppose $a$ and $b$ are homotopically non-trivial
simple loops 
in a  hyperbolic surface of injectivity radius $r$.
Then,

(a) (Thurston).  $ I(a,b) \leq \frac{4}{\pi r^2} l(a) l(b)$.

(b)  $l(D_a(b)) \leq I(a,b) l(a) + l(b)$.

(c) For each integer $n$, $\frac{\pi r^2  |n| I(a,b)}{ 4l(b)}
\leq l(D^n_a(b)) \leq |n| I(a,b) l(a) + l(b)$.

(d) If $|a \cap b| \geq 3$ or $|a \cap b| = 2$  so that the two points
of intersection have the same intersection signs, then there exists
a simple loop $c$ so that $l(c) \leq l(a) + l(b)$,
$ |D_c(b) \cap a| < |b \cap a|$ and $l(D_c(b)) \leq 2l(a) + l(b)$.

(e) There exists a sequence of simple loops $c_1, ..., c_k$ so that
$k \leq |a \cap b|$, $l(c_i) \leq (2i-1) l(a) + l(b)$ for each $i$ and
$D_{c_k} ... D_{c_1}(b)$ is either disjoint from $a$, or intersects $a$
at one point, or intersects $a$ at two points of different signs. 
\rm

\it Proof. \rm Part (a) is essentially in [FLP], pp.54, lemma 2.
We produce a slightly different proof so that the coefficient is
$\frac{4}{\pi r^2}$. Without loss of generality, we may assume that
both $a$ and $b$ are simple geodesics. 
Construct a flat torus as the metric product of two geodesics $a$ and $b$.
The area of the torus is $l(a) l(b)$. Each intersection point
of $a$ with $b$ gives a point $p$ in the torus. Now the flat
distance between any two of these points $p$'s  is at least the injectivity
radius $r$ (otherwise there would be Whitney discs for $a \cup b$).
Thus the flat disks of radius $r/2$ around these $p$'s
are pairwise disjoint. This shows that the sum of the areas of these
disks is at most $l(a)l(b)$ which is the Thurston's inequality.

To see part (b), we note that the Dehn-twisted loop $D_a(b)$ is
obtained by taking $I(a,b)$ many parallel copies of $a$ and resolving
all the intersection points between $b$ and the parallel copies (from
$a$ to $b$). Thus the inequality follows.

Part (c) follows from parts (a) and (b). Note that we have used
the fact that $I(D^n_a(b), b) = |n| I(a,b)$ (see for instance
[Lu] for a proof, or one also can check directly that there are
no Whitney disks for $D^n_a(b) \cup b$).

Part (d) is essentially in lemma 2 [Li]. Our minor observation is that
one can always choose  a positive Dehn-twist $D_c$ to achieve the
result.

We need to consider two cases.

\it Case 1. \rm There exist two intersection points $x$, $y$
$\in a \cap b$ adjacent along
in $a$ which have the same intersection signs (see figure 1).
Then
the curve $c$  as shown in figure 1 (with the right-hand
orientation on the surface)
satisfies  all conditions in the part (d). If the surface
is left-hand oriented, take $D_c(b)$ to be the  loop $c$.

\midspace{0.1cm}
\centerline{\epsfbox{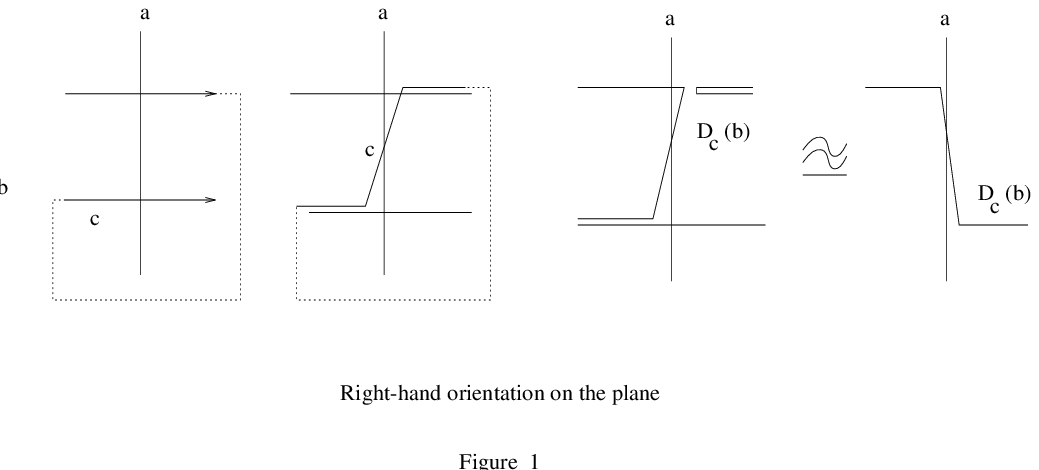}}
\midspace{0.1cm}

\midspace{0.1cm}
\centerline{\epsfbox{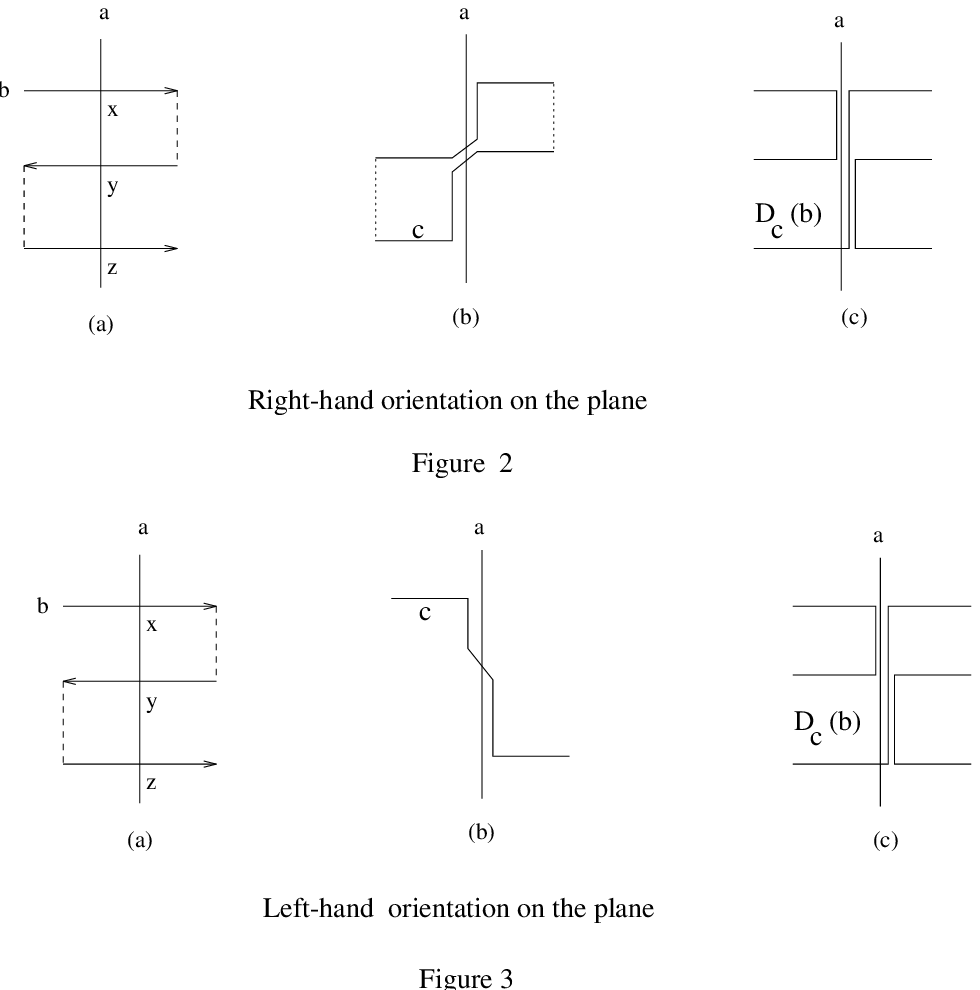}}
\midspace{0.1cm}

\it Case 2. \rm   Suppose any  pair of adjacent  intersection points
in $a \cap b$ has different intersection signs. Then $|a \cap b| \geq 3$.
Take three intersection
points $x, y, z \in a \cap b$ in so that $x$, $y$ and $y$, $z$ are adjacent
in $a$.
Their intersection signs alternate. 
Fix an orientation on $b$ so that
the arc from $x$ to $y$ in $b$ does not contain $z$ as shown in figure 2.
If  the surface $\Sigma$  is right-hand oriented as in figure 2,
take  $c$ as in figure 2(b). Then  $D_c(b)$ is shown in figure 2(c).
If the surface has the left-hand orientation, then take $c$ as
shown in figure 3(b). The loop $D_c(b)$ is shown in figure 3(c).
One checks easily that the simple loop $c$  satisfies the
all the conditions.

Part (e) follows from part (d) by induction on $|a \cap b|$.
$\square$

We shall also need the following well known lemma in order to deal
with disjoint loops and loops intersecting at one point.

{\bf Lemma 4.} \it Suppose $a$ and $b$ are two simple loops intersecting
transversely at one point. Then,

(a) $D_a D_b (a) \cong a$,

(b) $(D_a  D_b D_a)^2$ sends $a$ to $a$, $b$ to $b$ and
reverses the orientations  on both $a$ and $b$. 

\rm

See [Bir] and [Li] for a proof, or one can check it directly.
Note that  $(D_a  D_b D_a)^2$ is the hyper-elliptic involution
on the 1-holed torus containing both $a$ and $b$.

We first give a proof of theorem 2. The proof of theorem 1 follows
by making length estimate at each stage of the proof of theorem 2.

\it Proof of Theorem 2. \rm Let $f$ be an orientation preserving homeomorphism
of $\Sigma_{g,n}$ which is the identity map on the boundary.  We shall show that
there exists a composition  $p$ of positive Dehn-twists so that
for each $i$, $p f^{-1} |_{a_i} = id $. It follows that $pf^{-1}$ is
a product of Dehn-twists on $a_i$'s.  

We prove the theorem by induction on the norm $|\Sigma_{g,n}|
=3g-3+n$ of the surface (the norm is the complex dimension of
the Teichmuller space of complex structures with punctured ends
on the interior of the surface). The basic property of the norm
is that if $\Sigma'$ is an incompressible subsurface which is
not homotopic to $\Sigma_{g,n}$, then the norm of $\Sigma'$ is strictly
smaller than the norm of $\Sigma_{g,n}$. 
For simplicity, we assume that the Euler characteristic of the
surface is negative (though the proof below also works for the torus).

If the norm of a surface is zero, then the surface is the 3-holed
sphere. The theorem is known to hold in this case (see [De]).

If the norm of the surface is positive, we pick a non-boundary component,
say $a_1$, of the 3-holed sphere decomposition as follows. If the genus
of the surface $\Sigma_{g,n}$ is positive, $a_1$ is a non-separating loop. By
proposition 3(e) applied to $a= a_1$ and $b = f^{-1}(a_1)$, we find
a sequence of simple loops $c_1, ...., c_k$, $k \leq I(a,b)$ so that
$a'_1 = D_{c_k}...D_{c_1} f^{-1}(a_1)$ satisfies: either
$a_1' \cap a_1 = \emptyset$, or $|a'_1 \cap a_1| = I(a'_1, a_1) =1$,
or $|a'_1 \cap a_1| = 2_0$. There are two cases we need to consider:
(1) both $a_1$ and $a'_1$ are separating loops, and (2)
both of them are non-separating.

In the first case, by the choice of $a_1$, the genus of the surface is
zero. First of all $I(a_1, a'_1) = 1$ cannot occur due to homological
reason. Second, since the homeomorphism $D_{c_k} .. . D_{c_1} f^{-1}$
is the identity map on the non-empty boundary $\partial \Sigma_{0, n}$,
it follows that $I(a_1, a'_1) = 2$ is also impossible and
$a'_1$ is actually isotopic to $a_1$. After composing with an isotopy,
we may assume that $D_{c_k} ... D_{c_1} f^{-1} |_{a_1} $ is the identity
map. Now cut the surface open along $a_1$ to obtain two subsurfaces
of smaller norms. Each of these subsurfaces is stablized under 
$D_{c_k} ... D_{c_1} f^{-1}$. Thus induction hypothesis applies and
we conclude the proof in this case.

In the second case that both $a_1$ and $a'_1$ are non-separating, then
either $|a'_1 \cap a_1| = 1$, or there exists a third curve $c$ so that
$c$ transversely intersects  each of $a_1$ and $a'_1$ in one point.
By lemma 4(a), one of the product $h$ of positive Dehn-twists $D_{a'_1}
D_{a_1}$, or $D_cD_{a_1}D_{a'_1}D_c$ will send $a'_1$ to $a_1$.
If the homeomorphism $h D_{c_k} ... D_{c_1} f^{-1}$ sends $a_1$ to $a_1$
reversing the orientation, by lemma 4(b), we may use six more positive
Dehn-twists (on $a_1, a'_1$, or $c$, $a_1$) to correct the orientation.
Thus, we have constructed a composition of positive Dehn-twists
$D_{c_m}....D_{c_1} f^{-1}$ so that it is the identity map on $a_1$
and $m \leq I(a,b) + 10$. Now cut the surface open along $a_1$ and
use the induction hypothesis.  The result follows.
$\square$

We note that the proof fails if we do not choose $a_1$ to be a 
non-separating simple loop in the case the surface is closed of
positive genus.

Now we prove theorem 1 by making length estimate on each steps above.

\it Proof of Theorem 1. \rm
Let $f$ be an isometry of a hyperbolic closed surface $\Sigma_{g} = \Sigma_{g,0}$.

We begin with the following result which gives bound on the lengths of
$a_i$'s and $c  $ used in the proof of theorem 2.

{\bf Proposition 5.} \it Suppose $\Sigma_g$ is a hyperbolic surface of
injectivity radius $r$.

(a) (Bers) There exists a 3-holed sphere decomposition \{$a_1, ..., a_{3g-3}\}$
of the surface so that $l(a_i) \leq 26(g-1)$.

(b) If $a$ and $b$ are two non-separating simple geodesics in a compact
hyperbolic surface $\Sigma$ which is a totally geodesic subsurface
in $\Sigma_{g,n}$ so that
either $I(a, b) = 0$ or $| a \cap b|  = 2_0$, then there exists a simple
geodesic $c$ in $\Sigma$ so that $I(c,a) =I(c,b) =1$ and
$l(c) \leq \frac{8(g-1)r}{sinh r} + 8r$.

\rm

\it Proof. \rm See Buser [Bu], pp.123 for a proof of part (a).

To see part (b), we first note that there are simple loops $x$ so that
$I(x, a) =I(x,b) =1$ by the assumption on $a$ and $b$. Let $c$ be the
shortest simple loop in $\Sigma$ satisfying $I(c,a) = I(c,b) =1$.  We shall
estimate the length of $c$ as follows. Let $N =[\frac{l(c)}{2r}]$ be
the largest integer smaller than $\frac{l(c)}{2r}$. Let
$P_1 = a \cap c$, $P_2$, ..., $P_N$ be $N$ points in $c$ so that their
distances $d(P_i, P_{i+1}) = 2r$. Let  $B_i$ be the disc of radius $r$
centered at $P_i$ and $B_k$ be the ball containing $c \cap b$.
Then the shortest length  property of $c$ shows that the intersections of
the interior $int(B_i) \cap int( B_{j})$ is empty if $1 \leq i < j < k$ or
$k < i < j \leq N$. Thus the sum of the areas of the $N-2$ balls
$B_2, ..., B_{k-1}, B_{k+1}, ..., B_N$
is at most twice the area of the surface $\Sigma_{g,0}$.
This gives the estimate required.
$\square$

Fix a 3-holed sphere decomposition $\{a_1, ..., a_{3g-3}\}$ of
the hyperbolic surface so that $l(a_i) \leq 26(g-1)$. We may
assume that the loops $a_i$ are so labeled that
$a_1, a_2, ...., a_k$ are non-separating loops and the rest are separating.

We now show  that there exists a computerable constant $C'= C'(g,r)$
so that any orientation preserving isometry $f$ of the hyperbolic
surface $\Sigma_{g}$ is isotopic to a product $qp$ where $q$ is a 
product of   positive or negative Dehn-twists on $a_i$'s and 
$p$ is a product of
at most $C'(g,r)$ many positive Dehn-twists on curves of
lengths at most $C'(g,r)$.

We now rerun the  constructive  proof of 
theorem 2 by estimating the lengths of loops 
involved in the proof of theorem 2. To begin with, we take
$a = a_1$ and $b = f^{-1}(a_1)$ of lengths at most $26 g$. By
Thurston's inequality, their intersection number $I(a,b)
$ is at most 
$\frac{ 52^2 g^2}{\pi r^2}$.
By propositions 3(e), 5(b) and the proof of theorem 2,
we  produce a finite set of simple loops $\{c_1, ..., c_k\}$
so that $k \leq I(a,b) + 10$, the lengths of $c_i$ is bounded in
$g, r$ and
$f_1 =D_{c_k}....D_{c_1} f^{-1}$ is the identity  map on $a_1$.
Now we take $a = a_2$ and $b = f_1(a_2)$ and run the same
constructive proof as above in the totally geodesic subsurface
$\Sigma_{g-1, 2}$ obtained by cutting $\Sigma_{g}$ open along $a_1$.
In order for the proof to work, we need to see that the length
of $b$ is bounded. Indeed, proposition 2(b) gives the 
estimate of $l(b)$ in terms of $l(c_i)$,  $l(a_2)$, and
$g, r$ (here we estimate the intersection number $I(c_i, x)$
in terms of the lengths by Thuston's inequality).
Thus, we construct  a finite set of simple loops
$d_1, ..., d_m$ so that $m$ and $l(d_i)$ are bounded in $g,r$, 
$d_i \cap a_1 = \emptyset$, and $D_{d_m} ... D_{d_1} f_1^{-1}$ is the
identity map on $a_1 \cup a_2$. 
Inductively, we produce the required positive homeomorphism $p$.

We remark that if the injectivity radius $r$ is at least $log 2$, then
the number $C'(g,r)$ that we obtained is at least
$g^{g^ {g^{ .....^{g}}}}$ (there
are $3g-3$ many exponents) in magnitude.

As a consequence,  we obtain the following expression for the 
homeomorphism
$p^{-1} f = D^{n_1}_{a_1} ... D^{n_{3g-3}}_{a_{3g-3}}$. It remains to
show that the exponents $n_i$'s are bounded. To this end, for each
index $i$, we pick a geodesic loop $b_i$ which is disjoint from
all $a_j$'s for $j \neq i$ and $b_i$ intersects $a_i$ at one point
or two points of different signs. A simple calculation involving
right-angled hyperbolic hexagon shows that we can choose these
$b_i$ to have lengths at most $182(g-1) - log(r/4)$.  Thus the lengths of
curve $p^{-1} f(b_i)$ is bounded (in terms of $g$ and $r$). 
By proposition 3(d), the growth of the lengths of loops $D_{a_i}^n (b_i)$ is
linear in $|n|$ if $|n|$ is large. Thus we obtain an estimate on
the absolute value of the exponents $|n_i|$. This finishes the
proof.

\it Acknowledgment. \rm   The work is supported in part by the NSF.

\centerline{\bf References}

[Bi] Birman, J.: Mapping class groups of surfaces. In: Birman, J.,
Libgober, A. (eds.), Braids. Proceedings of a summer research
conference, Contemporary Math. Vol. 78, pp.13-44, Amer. Math. Soc. 1988

[Bu] Buser, P.: Geometry and spectra of compact Riemann surfaces. Progress
in Mathematics. Birkh\"auser, Boston, 1992

[De] Dehn, M.: Papers on group theory and topology. J. Stillwell (eds.).
 Springer-Verlag, Berlin-New York, 1987

[FLP] Fathi, A., Laudenbach, F., Poenaru, V.: Travaux de Thurston sur les
surfaces. Ast\'erisque {\bf 66-67}, Soci\'et\'e Math\'ematique de France, 1979

[Li] Lickorish, R.: A representation of oriented combinatorial 3-manifolds. Ann.
Math. {\bf 72} (1962), 531-540

[Lu] Luo, F.: Multiplication of simple loops on surfaces, preprint 1999

Department of Mathematics
 
Rutgers University

Piscataway, NJ 08903

\end